\documentclass[final,leqno,onefignum,onetabnum]{siamltex1213}

\newcommand{\tr}{\mathsf{T}}
\newcommand{\order}[1]{\mathbf{O}\left(#1\right)}

\usepackage[cmex10]{amsmath}
\usepackage{amssymb}
\usepackage{color}

\usepackage{graphicx}

\newtheorem{Remark}{Remark}[section]
\newtheorem{assumption}{Assumption}

\newcommand{\change}[1]{#1}

\title{Self-reflective model predictive control}

\author{Boris Houska$^{1,2}$ \and Dries Telen$^{3}$ \and Filip Logist$^{3}$ \and Jan Van Impe$^{3}$}

\begin{document}
\maketitle

\renewcommand{\thefootnote}{\arabic{footnote}}
\footnotetext[1]{Corresponding Author (\email{borish@shanghaitech.edu.cn}).}
\footnotetext[2]{School of Information Science and Technology, ShanghaiTech University, Shanghai, China.}
\footnotetext[3]{Chemical Engineering Department, KU Leuven, BioTeC+ \& OPTEC, Gebroeders De Smetstraat, 9000 Ghent, Belgium.}

\begin{abstract}
This paper proposes a novel control scheme, named self-reflective model predictive control, which takes its own limitations in the presence of process noise and measurement errors into account. In contrast to existing output-feedback MPC and persistently exciting MPC controllers, the proposed self-reflective MPC controller does not only propagate a matrix-valued state forward in time in order to predict the variance of future state-estimates, but it also propagates a matrix-valued adjoint state backward in time. This adjoint state is used by the controller to compute and minimize a second order approximation of its own expected loss of control performance in the presence of random process noise and inexact state estimates. The properties of the proposed controller are illustrated with a small but non-trivial case study.
\end{abstract}




\pagestyle{myheadings}
\thispagestyle{plain}

\section{Introduction}

Historically, online optimization based receding horizon controllers have been developed by separating the overall control problem into two, namely, a state estimation and a state feedback control problem~\cite{Allgower2000,Rawlings2009}. Nowadays, there exists a mature theoretical foundation as well as software for both moving horizon estimation (MHE)~\cite{Diehl2009,Rao2001}, realizing optimized state estimation, as well as model predictive control (MPC)~\cite{Bemporad2002,Biegler1991,Chen1998,Houska2011}, realizing optimized state feedback control.

Unfortunately, if the dynamic system is not linear or if there are state- or control constraints presents, a separation of the overall control problem into an estimation and a feedback phase can be sub-optimal. For example, if a nonlinear system is not observable at its desired operation point, it might be necessary to excite the system from time to time on purpose in order to be able to measure the states. One way to address this problem is by using a so-called output-feedback model predictive controller, as suggested in~\cite{Findeisen2003,Mayne2006,Mayne2009}. During the last decade, there has been significant progress in the field of output-feedback MPC. In particular, in~\cite{Goulart2007} a way of approaching the output-feedback MPC problem for linear systems in the presence of convex state- and control constraints has been developed. The corresponding method uses convex optimization techniques in order to design a feedback law that is affine in the disturbance sequence by solving an online optimization problem whose complexity scales quadratically with respect to the prediction horizon. Other approaches for conservatively solving the output-feedback problem based on set propagation techniques can be found in~\cite{Bemporad2000,Chisci2002,Kurzhanski2004,Lee2001,Yan2005}.

Notice that the output-feedback model predictive control problem is related to the so-called \textit{dual control} problem~\cite{Feldbaum1961}. Dual control tries to find an optimal control input under the assumption that learning is possible while the model is unknown. A complete overview about the historical developments related to the dual control problem is beyond the scope of this article, but there exists a large number of recommendable overview articles~\cite{Filatov2000,Filatov2004,Lee2009,Wittenmark1985,Wittenmark1995}, where reviews and comments about the history and importance of the dual control problem can be found. An early article about the limitations of predictive control from an adaptive control perspective can be found in~\cite{Bitmead1991}. During the last decade, a number of articles have appeared, which all try in one or the other way to include dual control aspects in MPC, as reviewed below. Here, one way to predict future state- and parameter estimation errors is to use an extended Kalman filter.
This idea has been introduced by Hovd and Bitmead~\cite{Hovd2005}, who propose to augment standard MPC by an additional matrix-valued state that predicts the variance of future state estimates. This matrix-valued state is propagated forward in time by a Riccati differential equation, which computes the variance state of an extended Kalman filter~\cite{Telen2013}. It can be shown that the corresponding controller excites the system for improved state estimates. A related approach toward dual MPC has been proposed in~\cite{Heirung2012,Heirung2015}, where the MPC objective is augmented by a term that---similar to optimal experiment design~\cite{Fisher1935,Pukelsheim1993,Telen2014}---penalizes the predicted parameter error variance. In~\cite{Heirung2015a} the same authors develop a way to predict and optimize future parameter estimation errors for single-input-single-output finite impulse response systems based on non-convex quadratically constrained quadratic programming formulations. A different way of taking predicted parameter estimation errors into account has been suggested in~\cite{Kim2008}, where an MPC scheme is augmented by an adaptive parameter estimation algorithm. Yet another related control problem formulation is proposed in~\cite{Larrson2013}, where it is suggested to add an experiment design constraint to the MPC formulation. Lucia and Paulen~\cite{Lucia2014} suggest to combine guaranteed parameter estimation and multi-stage MPC techniques in order to reduce the uncertainty about future estimation errors. Moreover, in~\cite{Forgione2015} and \cite{Larsson2015} various extensions of application oriented experiment design are developed in order to generate excitation signals in the context of receding horizon control. Another recent article~\cite{Mesbah2015} on online experiment design suggests to use a modified experiment design framework for generating ``least costly'' excitations, which try to avoid perturbation of the nominal plant operation as much as possible. Other approaches for persistently exciting closed loop MPC can be found in~\cite{Hernandez2015,Marafioti2014,Shouche2002,Zacekova2013}.

\subsection{Contribution}

After introducing the mathematical notation and background in Section~\ref{sec::prelim}, the paper starts with Section~\ref{sec::motivating_example} by outlining an example, which illustrates why MPC may fail to stabilize nonlinear control systems if the MPC controller does not take the properties of the state estimator into account. Sections~\ref{sec::mpc} and~\ref{sec::mhe} review existing concepts from the field of MPC, dynamic programming, and state estimators based on extended Kalman filters and MHE. The main theoretical contribution of this paper is presented in Section~\ref{sec::loss}, which proposes a novel way of analyzing the loss of optimality that is associated with implementing certainty equivalent MPC controllers in the presence of process disturbances and measurement noise. This analysis is related to a recently proposed method for formulating economic optimal experiment design problems~\cite{Houska2015}. However, \cite{Houska2015} analyzes a static optimization problem while this paper is concerned with the loss of optimality of a dynamic two-player-multi-stage game, as the MPC controller is optimizing the control input whenever new state estimates become available. A major contribution in this context is presented in Theorem~\ref{thm::expansion}, which proposes a procedure to compute second order moment expansions of the expected loss of optimality of certainty equivalent MPC in the presence of random process noise as well as random measurement errors. The corresponding method scales linearly with the prediction horizon of the MPC controller despite the fact that the future state estimates are correlated in time. In Section~\ref{sec::rhc} this second order moment approximation of the expected sub-optimality of certainty equivalent MPC is used to design a novel ``self-reflective'' model predictive controller, which optimizes not only its nominal performance but also includes a term that approximates its own expected sub-optimality in the presence of process noise and measurement errors. This leads---similar to the above reviewed dual MPC schemes---to a controller that persistently excites the states in order to reduce the variance of future state estimates. However, the above reviewed existing dual MPC schemes are all based on additional matrix-valued hyperstates that propagate forward in time in order to predict the accuracy of future state estimates. The proposed receding horizon controller also propagates states forward in time, but, in contrast to existing dual MPC schemes, it is based on additional matrix-valued adjoint states, which propagate backward in time. By coupling the information of these forward and adjoint states the control scheme obtains the ability to be ``self-reflective'', in the sense that it is capable of taking its own limitations into account. The novelty and differences of this approach compared to the numerous existing articles on persistently exciting MPC is further discussed toward the end of Section~\ref{sec::rhc}.  Section~\ref{sec::numerics} presents a numerical case study, which illustrates the properties of self-reflective MPC. Section~\ref{sec::conclusions} concludes the paper.

\subsection{Notation}
\label{sec::notation}
For any matrix $M \in \mathbb R^{n \times n}$, the notation $M \succ 0$ indicates that $M$ is positive definite. The notation $M^\dagger$ is used to denote the \change{pseudo-inverse} of a matrix $M$. If $\varphi: \mathbb R^{n} \to \mathbb R^{m}$ is an integrable function and $x$ a random variable with integrable probability distribution $\rho: X \to \mathbb R$ and support $X \subseteq \mathbb R^n$, the expected function value of $\varphi$ is denoted by
$$\mathbb E_x \{ \varphi(x) \} \; = \; \int_{X} \, \varphi(x) \rho(x) \, \mathrm{d}x \; .$$
Notice that the index is occasionally omitted, i.e., we simply write $\mathbb E \{ \varphi(x) \}$, if it is clear from the context that $x$ is the random variable with respect to which the expected value is computed. If the random variable has a block structure,
$$x = \change{\left[ x_1^\tr \, x_2^\tr \right]}^\tr \quad \text{with} \quad X_1(x_2) = \left\{ x_1 \in \mathbb R^{n_1} \mid \change{ \left[ x_1^\tr \, x_2^\tr \right]^\tr } \in X \right\}$$
and $0 \leq n_1 \leq n$, the conditional expectation over $x_1$ is denoted by
$$\mathbb E_{x_1} \{ \varphi(x) \} \; = \; \int_{X_1(x_2)} \, \varphi(x) \rho(x) \, \mathrm{d}x_1 \; .$$
Notice that $\mathbb E_{x_1} \{ \varphi(x) \}$ is a function of $x_2$. In particular, this notation is such that
\[
\mathbb E_{x_1} \mathbb E_{x_2} \{ \varphi(x) \} = \mathbb E_{x_2} \mathbb E_{x_1} \{ \varphi(x) \} = \mathbb E_{x} \{ \varphi(x) \} \; .
\]
In addition, the notation $\Vert \cdot \Vert_2$ denotes either the Euclidean norm, if the argument is a vector, or the spectral norm, if the argument is a quadratic matrix. Moreover, we use the upper index notation
$$x^{[k]} = \change{ \left[ x_0 \, x_1 \, x_2 \, \ldots \, x_k \right] }$$
to denote the sub-matrix that consists of the first $k+1$ columns of a given matrix $x = \change{ \left[ x_0 \, x_1 \, \ldots \, x_N \right] }$ with $x_0,x_1,\ldots,x_N \in \mathbb R^{n}$ being column vectors. This upper index notation is well-defined for all integers $0 \leq k \leq N$. Additionally, we allow the upper index ``$^{[-1]}$'', which should be read as ``none''. For example, if we write an expression like
$$\varphi_k(x^{[k-1]}) \; , \; k \in \{ 0, \ldots, N \} \; , $$
this means that the $k$-th element of the function sequence $\varphi_0, \varphi_1, \ldots, \varphi_N$ depends on the first $k$ columns of $x$ if $k \geq 1$, but the function $\varphi_0$ depends on \textit{none} of the elements of $x$. Similarly, the lower index notation
$$x_{[k]} = \change{ \left[ x_k \, x_{k+1} \, \ldots \, x_N \right] }$$
is used to denote the sub-matrix that consists of all columns of $x$ with index larger than or equal to $k$. Notice that throughout this paper $x$ denotes a state sequence of a discrete-time system, whose index typically runs from $0$ to $N$. The associated control and disturbance sequences are denoted by $u$ and $w$ and have one element less. For example,
$$u_{[k]} = \change{ \left[ u_k \, u_{k+1} \, \ldots \, u_{N-1} \right] } \quad \text{and} \quad w_{[k]} = \change{ \left[ w_k \, w_{k+1} \, \ldots \, w_{N-1} \right] }$$
denotes the sequences of controls and disturbances from $k$ to $N-1$. Analogous to the upper index notation, we allow an ``index overflow'', i.e., the lower index ``$_{[N]}$'' should be read as ``none''.

\section{Preliminaries}
\label{sec::prelim}

This paper concerns the design of control laws for nonlinear \change{discrete-time} systems of the form
\begin{align}
\label{eq::nonlinearSystem}
\begin{array}{rcl}
x_{k+1} &=& f(x_k,u_k) + w_k \\
\eta_k &=& h(x_k) + v_k \; .
\end{array}
\end{align}
Here, $x_k \in \mathbb R^{n_x}$ denotes the state, $u_k \in \mathbb R^{n_u}$ the control input, and $w_k \in \mathbb R^{n_x}$ the process noise at time $k \in \mathbb Z$. The variable $\eta_k \in \mathbb R^{n_h}$ denotes the output measurement and $v_k \in \mathbb R^{n_h}$ the measurement error. The functions
$$f: \mathbb R^{n_x} \times \mathbb R^{n_u} \to \mathbb R^{n_x \times n_x} \quad \text{and} \quad h: \mathbb R^{n_x} \to \mathbb R^{n_h}$$
denote the right-hand side function of the discrete-time dynamic system as well as the measurement function. Moreover, the stage- and terminal cost are denoted by
$$l: \mathbb R^{n_x} \times \mathbb R^{n_u} \to \mathbb R \quad \text{and} \quad m: \mathbb R^{n_x} \to \mathbb R \; .$$
Throughout this paper, the following blanket assumption is used.

\bigskip
\begin{assumption}
\label{ass::differentiable}
The functions $f,h,l$, and $m$ are three times continuously differentiable with respect to all arguments.
\end{assumption}

\bigskip
\noindent
At many places in the paper, shorthands are used to denote derivatives of the functions $f$, $h$, and $l$. For the first order derivatives the following shorthands are used
\[
\begin{array}{rclrclrcl}
A_k &=& \dfrac{\partial f(x_k,u_k)}{\partial x}, \quad & B_k &=& \dfrac{\partial f(x_k,u_k)}{\partial u}, \quad & C_k &=& \dfrac{\partial h(x_k)}{\partial x}, \\[0.3cm]
q_k &=& \nabla_x l(x_k,u_k) \; .
\end{array}
\]
Notice that this notation will only be used, if it is clear from the context at which point the derivatives of the functions $f$, $h$, and $l$ are evaluated. Similarly, the associated second order derivatives of these functions are denoted by
\[
\begin{array}{rclrclrcl}
K_k &=& \nabla_{xx}^2 f(x_k,u_k), & L_k &=& \nabla_{ux}^2 f(x_k,u_k), \quad & M_k &=& \nabla_{uu}^2 f(x_k,u_k), \\[0.2cm]
Q_k &=& \nabla_{xx}^2 l(x_k,u_k), & R_k &=& \nabla_{uu}^2 l(x_k,u_k), \quad & S_k &=& \nabla_{xu} l(x_k,u_k) \; .
\end{array}
\]
All technical derivations are based on the following assumptions about the process noise and measurement errors.

\bigskip
\begin{assumption}
\label{ass::variances}
The process noise variables $w_k$ and the measurement errors $v_k$, $k \in \mathbb N$, are independent random variables, whose expectation values and variances $W \succeq 0$ and $V \succeq 0$ are given,
\[
\begin{array}{rclrcl}
\mathbb E \{ w_k \} &=& 0 \; , \quad & \mathbb E \{ w_k w_k^\tr \} &=& W \\[0.1cm]
\mathbb E \{ v_k \} &=& 0 \; , \quad & \mathbb E \{ v_k v_k^\tr \} &=& V \; .
\end{array}
\]
\end{assumption}

\bigskip

\begin{assumption}
\label{ass::boundedSupport}
The process noise variables $w_k$ and the measurement errors $v_k$ have bounded support such that $\Vert v_k \Vert_2 \leq \gamma$ \change{and} $\Vert w_k \Vert_2 \leq \gamma$ for a given radius $\gamma > 0$ such that $$\Vert V \Vert_2 \leq \gamma^2 \quad \change{\text{and}} \quad \Vert W \Vert_2 \leq \gamma^2 \; .$$
\end{assumption}

\bigskip
\noindent
Notice that Assumptions~\ref{ass::variances} and ~\ref{ass::boundedSupport} together allow us to compute second order moment expansions of nonlinear functions. For example, if $\varphi: \mathbb R^{n_x} \to \mathbb R$ is a three times continuously differentiable function of $w_k$ with $\varphi(0) = 0$, we have~\cite{Oehlert1992}
\begin{eqnarray}
\mathbb E \{ \varphi(w_k) \} &=& \mathbb E \left\{ \varphi(0) + \nabla \varphi(0) w_k + \frac{1}{2} w_k^\tr ( \nabla^2 \varphi(0) ) w_k + \order{\gamma^3} \right\} \notag \\[0.1cm]
&=& \mathbb E \left\{ \frac{1}{2} \mathrm{Tr} \left( \nabla^2 \varphi(0) w_k w_k^\tr \right) \right\} + \order{\gamma^3} \; . \notag
\end{eqnarray}
This equation can also be written in the form
\begin{eqnarray}
\label{eq::expansion1}
\mathbb E \{ \varphi(w_k) \} &=& \frac{1}{2} \mathrm{Tr}( \nabla^2 \varphi(w') W ) + \order{\gamma^3} \; ,
\end{eqnarray}
which holds independently of the choice of the point $w'$, at which the derivative $\nabla^2 \varphi$ is evaluated, as long as $\Vert w' \Vert_2 = \order{\gamma}$. Similarly, the equation
\begin{eqnarray}
\mathbb E \{ \varphi(w_k)^2 \} &=& \mathbb E \left\{ \left( \varphi(0) + \nabla \varphi(0) w_k + \frac{1}{2} w_k^\tr ( \nabla^2 \varphi(0) ) w_k + \order{\gamma^3} \right)^2  \right\} \notag \\[0.1cm]
&=& \nabla \varphi(0)^\tr W \nabla \varphi(0) + \order{\gamma^3} \notag
\end{eqnarray}
with $W = \mathbb E \{ w_k w_k^\tr \}$ implies that we have
\begin{eqnarray}
\label{eq::expansion2}
\mathbb E \{ \varphi(w_k)^2 \} &=& \nabla \varphi(w')^\tr W \nabla \varphi(w') + \order{\gamma^3} \; .
\end{eqnarray}
Also this result holds independently of the choice of the point $w'$ as long as we ensure that $\Vert w' \Vert_2 = \order{\gamma}$. Similar expansion results hold for three times continuously differentiable vector-valued functions. For example, if $\tilde \varphi: \mathbb R^{n_x} \to \mathbb R^{n_\varphi}$ is a three times continuously differentiable function, we have
\begin{eqnarray}
\mathrm{Var}(\tilde \varphi(w_k) ) &=& \mathbb E \left\{ \left( \tilde \varphi(w_k) - \mathbb E \{ \tilde \varphi(w_k) \} \right)\left( \tilde \varphi(w_k) - \mathbb E \{ \tilde \varphi(w_k) \} \right)^\tr   \right\} \notag \\[0.16cm]
\label{eq::expansion3}
&=& \nabla \tilde \varphi(w')^\tr W \nabla \tilde \varphi(w') + \order{\gamma^3}
\end{eqnarray}
for any $w'$ with $\Vert w' \Vert_2 = \order{\gamma}$. We will use equations of the form~\eqref{eq::expansion1},~\eqref{eq::expansion2}, and~\eqref{eq::expansion3} in this or very similar versions at many places in this paper (without always mentioning this explicitly).

\bigskip
\begin{Remark}
In practice, random variables are often modelled by Gaussian probability distributions. Unfortunately, Assumption~\ref{ass::boundedSupport} is violated in this case, as Gaussian distributions do not have a bounded support. It is possible to rescue\change{~\eqref{eq::expansion1} and~\eqref{eq::expansion2}} for Gaussian distributions with small variance, $\Vert W \Vert_2 \leq \gamma$, if we impose more restrictive assumptions on the global properties of the function $\varphi$. For example, if $\varphi$ has uniformly bounded third derivatives on $\mathbb R^{n_x}$, \change{~\eqref{eq::expansion1} and~\eqref{eq::expansion2} hold}. However, for general nonlinear functions $\varphi$ the equations~\eqref{eq::expansion1} and~\eqref{eq::expansion2} are wrong. Therefore, this paper suggests to model the random variables with probability distributions that closely approximate Gaussian distributions, but have bounded support. One might argue, that this is a realistic modeling assumptions when dealing with nonlinear systems, as noise terms \mbox{are---due} to physical \mbox{limitations---in} practice almost always bounded anyhow, which justifies Assumptions~\ref{ass::boundedSupport}.\hfill$\diamond$
\end{Remark}

\section{Motivating Example}
\label{sec::motivating_example}
Let us consider a nonlinear \change{discrete-time} system of the form~\eqref{eq::nonlinearSystem} with dimensions $n_x = n_u = 2$ as well as $n_h = 1$, given by
\begin{align}
\label{eq::f}
\left\{
\begin{array}{rcl}
f(x_k,u_k) &=& F(u_k) x_k + \delta u_k \\[0.1cm]
h(x_k) &=& C x_k
\end{array}
\right\}
\qquad \text{with} \qquad F(u_k) = 
I + \delta \left( \begin{array}{cc}
1 & 0 \\
C u_k & -1
\end{array}
\right)
\end{align}
and $C = (0,1)$. Here, $\delta \ll 1$ is a small time step. This is a nonlinear system as the system matrix $F(u_k)$ depends on the control input $u_k$. The stage cost is in this example assumed to be a tracking term of the form
\begin{align}
\label{eq::l}
l(x_k,u_k) = \delta \left( x_k^\tr Q x_k + u_k^\tr R u_k \right) \qquad \text{with} \quad R = I \quad \text{and} \quad Q = \left(
\begin{array}{cc}
1 & 0 \\
0 & 0
\end{array}
\right)\; .
\end{align}
It can be checked that this system is unstable in open-loop mode. Since only the second state component, $C x_k$, can be measured, the only way to gather information about the first component is to choose $u \neq 0$ such that the first state component has an influence on the propagation of the second state component via the term $C u_k$ in the lower left corner of the system matrix $F(u_k)$. If we choose $u$ such that $C u_k = 0$ for all $k$, the first component of the state vector $x$ cannot be observed. Unfortunately, a standard tracking MPC controller, whose aim is to minimize the stage cost $l$ under the assumption that full-state measurements are available, will choose $C u_k = 0$ for all $k$. Thus, the above nonlinear system cannot be stabilized if the MPC controller does not cooperate with the state estimator. For example, a naive cascade consisting of an extended Kalman filter (or MHE) and MPC leads to unstable closed-loop trajectories. The goal of this paper is to develop a receding horizon controller, which resolves this problem while maintaining economic performance.

\section{Model Predictive Control}
\label{sec::mpc}

This section briefly reviews the main idea of model predictive control (MPC) and its associated cost-to-go function~\cite{Rawlings2009}. For this aim, we introduce for all $i \in \{ 0,1, \ldots, N \}$ the notation
\begin{align}
\label{eq::mpc}
\begin{array}{rccl}
J_i(x_i,w_{[i]}) &=& \underset{x,u}{\text{min}} \; & \sum_{k=i}^{N-1} l(x_k,u_k) + m(x_N) \\[0.3cm]
& & \text{s.t.} \; &
x_{k+1} = f(x_k,u_k) + w_k \; , \; \; \change{\text{for all} \; \; } k \in \{ i, \ldots, N-1 \} \\
\end{array}
\end{align}
to denote the parametric optimal control problem that is associated with the nonlinear system $f$, stage cost $l$ and terminal cost $m$. The function $J_i: \mathbb R^{n_x} \times \mathbb R^{n_x \times (N-i)} \to \mathbb R$ is called the cost-to-go function that is associated with the time horizon $\{ i, i+1, \ldots, N \}$. Here, $w_{[i]} = \change{ \left[ w_i \, \ldots \, w_{N-1} \right] }$ denotes the last $(N-i)$ elements of the disturbance sequence as explained in Section~\ref{sec::notation}. Notice that the last cost-to-go function, $J_N = m$, is equal to the Mayer term $m$. The standard implementation of an MPC loop can be outlined as follows.
\bigskip
\begin{enumerate}
\item Wait for an estimate $y_0$ of the current state.
\item Solve Problem~\eqref{eq::mpc} for $i=0$, $x_0 = y_0$, and $w = \change{ \left[ w_0 \, w_1 \, \ldots \, w_{N-1} \right] } = 0$.
\item Send the first element $u_0$ of the optimal control sequence to the real process.
\item Shift the time index and continue with Step 1.
\end{enumerate}
\bigskip
Notice that MPC is a certainty equivalent controller. This means that the future control input is optimized as if there was no uncertainty, $w = 0$, and as if the state estimate $x_0 = y_0$ was accurate. Comments on how to choose the terminal cost $m$ and how to ensure stability of MPC can be found in~\cite{Chen1998,Mayne2000,Rawlings2009}.

\bigskip
\begin{Remark}
In the following discussion control and/or state constraints are not analyzed explicitly, as this article focuses on how to design controllers for the case that the function $f$, $l$, and $m$ are nonlinear.
The considerations in the following sections can, however, be generalized for the case that control constraints are present, as explained at the end of Section~\ref{sec::rhc}. State constraints are beyond the scope of this paper.\hfill$\diamond$
\end{Remark}

\bigskip
\noindent
In the following, we outline the main idea of \textit{dynamic programming}~\cite{Bellman1957,Bertsekas2012}. For this aim, we introduce the auxiliary functions
\begin{align}
\label{eq::Gdef}
G_k(x_k,u_k,w_{[k]}) = J_{k+1}(f(x_k,u_k)+w_k,w_{[k+1]}) + l(x_k,u_k) \; ,
\end{align}
which are defined for all $x_k \in \mathbb R^{n_x}$, $u_k \in \mathbb R^{n_u}$ and all $w_{[k]} \in \mathbb R^{n_x \times (N-k)}$, and all $k \in \{ 0,1, \ldots, N-1\}$. The cost-to-go function sequence $J_0,J_1,\ldots,J_N$ satisfies a functional recursion, which is known under the name dynamic programming. It starts with $J_N(\cdot) = m(\cdot)$ and iterates backwards
\begin{align}
\label{eq::dynProg}
J_{k}(\cdot,w_{[k]}) = \min_{u_k} \; G_k(\cdot,u_k,w_{[k]})
\end{align}
for all $k \in \{ 0, 1, \ldots, N-1 \}$. Notice that this functional recursion is useful for analyzing the properties of MPC~\cite{Rawlings2009}, although numerical implementations of MPC typically solve problem~\eqref{eq::mpc} directly instead of using dynamic programming.
Assumption~\ref{ass::differentiable} implies that the function $J_N = m$ as well as the function $G_{N-1}$ are three times continuously differentiable. The following proposition introduces a regularity condition under which this property is inherited by all cost-to-go functions through the dynamic programming recursion.

\bigskip
\begin{proposition}
\label{prop::diff}
Let Assumption~\ref{ass::differentiable} be satisfied and let $k \in \{ 0, \ldots, N-1 \}$ be given. If the function $J_{k+1}$ is locally three times continuously differentiable in a neighborhood of $(x_k,0)$, then the function $G_k$ is locally three times continuously differentiable in all arguments. Moreover, if the minimizer $\mu_k(x_k)$ of problem~\eqref{eq::dynProg} is unique and satisfies the first order necessary and second order sufficient optimality conditions,
\begin{align}
\label{eq::stationarity}
& \nabla_u G_k(x_k,\mu_k(x_k),0) = 0 \\[0.16cm]
\label{eq::regularity}
\text{and} \quad & G_{k}^{uu} = \nabla_{uu}^2 G_k(x_k,\mu_k(x_k),0) \succ 0
\end{align}
for all $x_k$, then the functions $J_{k}$ is locally three times continuously differentiable in a neighborhood of $(x_k,0)$.
\end{proposition}

\bigskip
\noindent
\proof
The first statement follows trivially from the definition of $G_k$, as the composition and sum of three times continuously differentiable functions remains three times continuously differentiable. The second statement is a well-known result from the field of parametric nonlinear programming~\cite{Nocedal2006}, which can be established by applying the implicit function theorem to the first order optimality condition~\eqref{eq::stationarity}.\hfill$\diamond$

\bigskip
\noindent
Notice that Proposition~\ref{prop::diff} can be applied recursively, in order to show that the cost-to-go functions $J_k$ as well as the functions $G_k$ are---under the regularity assumptions from Proposition~\ref{prop::diff}---for all $k$ locally three times differentiable. In this context, the word ``locally'' means that the functions $J_k$ and $G_k$ are three times continuously differentiable for small disturbance sequences $w_{[k]}$. This is sufficient, as we are in this paper interested in analyzing the influence of small disturbances, i.e., for the case that Assumption~\ref{ass::boundedSupport} is satisfied for a sufficiently small $\gamma > 0$.

\bigskip
\begin{assumption}
\label{as::regularity}
We assume that the regularity conditions from Proposition~\ref{prop::diff} are satisfied recursively for all $k \in \{ N-1,N-2,\ldots,0 \}$ such that all cost-to-go functions $J_k$ are locally three-times continuously differentiable for all $k \in \{ 0, \ldots, N \}$.
\end{assumption}

\bigskip
\noindent
In the following, we introduce a shorthand notation for the first and second order derivatives of the cost-to-go functions $J_k$,
\[
p_k = \nabla_x J_k(x_k,0) \quad \text{and} \quad P_k = \nabla_x^2 J_k(x_k,0) \; .
\]
Notice that the sequence $\Omega_k = (p_k,P_k)$ satisfies an algebraic Riccati recursion~\cite{Bittanti1991} of the form
\begin{eqnarray}
\label{eq::forwardpP}
\forall k \in \{ 0, \ldots, N-1\}, \quad  \Omega_{k} = \mathcal F_{\Omega}(x_k,u_k,\Omega_{k+1}) \quad \text{and} \quad \Omega_N = \mathcal M(x_N) \; .
\end{eqnarray}
Here, the function
$$\mathcal M(x_N) \; = \; \change{ \left[ \nabla_x J_N(x_N) \, \nabla_x^2 J_N(x_N) \right] } \; = \; \change{ \left[ \nabla_x m(x_N) \, \nabla_x^2 m(x_N) \right] } \in \mathbb R^{(n_x+1) \times n_x}$$
is for all $x_N \in \mathbb R^{n_x}$ given by the associated first and second order derivative of the Mayer term $m$.
The expression for the right-hand side function $\mathcal F_{\Omega}$ can be worked out explicitly by differentiating the recursion~\eqref{eq::dynProg} twice. This yields
\begin{align}
\label{eq::backwardpP}
\mathcal F_{\Omega}(x_k,u_k,\Omega_{k+1}) = \change{ \left[ \; A_k^\tr p_{k+1} + q_k \; \; 
G_{k}^{xx} - G_{k}^{xu} \left( G_{k}^{uu} \right)^{-1} G_{k}^{ux} \; \right] } \in \mathbb R^{(n_x+1) \times n_x}.
\end{align}
Here, $G_k^{xx}, G_k^{ux} = \left( G_k^{xu} \right)^\tr$, and $G_k^{uu}$ are shorthands for the second order derivatives of the function $G_k$ with respect to $x$ and $u$,
\begin{align}
\label{eq::G}
\begin{array}{rcl}
G_{k}^{xx} &=& Q_k + A_k^\tr P_{k+1} A_k + K_k \cdot p_{k+1} \\
G_{k}^{ux} &=& S_k + B_k^\tr P_{k+1} A_k + L_k \cdot p_{k+1} \\
G_{k}^{uu} &=& R_k + B_k^\tr P_{k+1} B_k + M_k \cdot p_{k+1} \; .
\end{array}
\end{align}
Recall that the shorthands $A_k$, $B_k$, $K_k$, $L_k$, $M_k$, $Q_k$, $S_k$, and $R_k$ denote first and second order derivatives of the functions $f $ and $l$ with respect to states and controls as defined in Section~\ref{sec::prelim}. The first two arguments of the function $\mathcal F_{\Omega}(x_k,u_k,\Omega_{k+1})$ on the left hand side of~\eqref{eq::backwardpP} indicate that all derivatives are evaluated at the point $(x_k,u_k)$, although this dependence is not highlighted explicitly in the right-hand side of~\eqref{eq::backwardpP}. Finally, a second order expansion of the cost-to-go functions is given by
\begin{eqnarray}
\label{eq::secondOrderJ}
J_k(x,0) &=& J_k(x_k,0) + p_k^\tr (x-x_k) + \frac{1}{2} (x-x_k)^\tr P_k (x-x_k) + \order{\Vert x-x_k \Vert^{3}}
\end{eqnarray}
for all $k \in \{ 0, 1, \ldots, N \}$.

\bigskip
\begin{Remark}
If $f$ is \change{affine} in $(x,u)$ while $l$ and $m$ are quadratic forms, we have $K_k = 0$, $L_k = 0$, and $M_k  = 0$,  the functions $J_k$ are quadratic forms in $x$, and the second order expansion~\eqref{eq::secondOrderJ} is exact. In this case, the algebraic Riccati recursion~\eqref{eq::forwardpP} for the sequence $\Omega_k = (p_k,P_k)$ has the form
\begin{eqnarray}
p_{k} &=& A_k^\tr p_{k+1} + q_k \notag \\
P_k &=& Q_k + A_k^\tr P_{k+1} A_k - \left( S_k + B_k^\tr P_{k+1} A_k \right)^\tr \left( R_k + B_k^\tr P_{k+1} B_k \right)^{-1} \left( S_k + B_k^\tr P_{k+1} A_k \right) \notag
\end{eqnarray}
and the corresponding MPC controller is equivalent to LQR control~\cite{Bittanti1991}.\hfill$\diamond$
\end{Remark}

\section{State Estimation}
\label{sec::mhe}
This section briefly reviews the extended Kalman filter (EKF)~\cite{Stengel1994} and its relation to moving horizon estimation (MHE)~\cite{Rawlings2009}. Recall that $\eta_k = h(x_k)+v_k$ denotes the measurements at time~$k$. The EKF proceeds by maintaining a state estimate $y_k \in \mathbb R^{n_x}$ and a variance $\Sigma_k \in \mathbb R^{n_x \times n_x}$ by performing updates of the form
\begin{eqnarray}
\label{eq::ekf1}
y_{k+1} &=& \mathcal F_y( x_{k}, u_{k}, y_k, \eta_k ) \; ,\\[0.16cm]
\label{eq::ekf2}
\Sigma_{k+1} &=& \mathcal F_{\Sigma}( x_{k}, u_{k}, \Sigma_k )
\end{eqnarray}
where $(x_{k},u_{k})$ denotes the point at which the nonlinear functions $f$ and $h$ are linearized. In detail, the right-hand side functions $\mathcal F_{y}$ and $\mathcal F_{\Sigma}$ of the EKF are given by the explicit expressions~\cite{Stengel1994}
\begin{eqnarray}
\mathcal F_y( x_{k}, u_{k}, y_k, \eta_k ) &=&  f \left( y_k + \Sigma_k C_k^\tr \left( C_k \Sigma_k C_{k}^\tr + V \right)^\dagger \left( \eta_k - h(y_k) \right) , u_k \right) \; , \notag \\[0.16cm]
\mathcal F_{\Sigma}( x_k, u_k, \Sigma_k ) &=& A_k \left( \Sigma_k - \Sigma_k C_k^\tr \left( C_k \Sigma_k C_k^\tr + V \right)^\dagger C_k \Sigma_k \right) A_k^\tr + W \; . \notag 
\end{eqnarray}
The most common variant of the EKF evaluates the derivatives
$$A_k = \frac{\partial f(x_k,u_k)}{\partial x} \quad \text{and} \quad C_k = \frac{\partial h(x_k)}{\partial x}$$
at the current state estimate $x_k = y_k$. However, in the following, we keep our notation general remarking that the first two arguments of the functions $\mathcal F_y$ and $\mathcal F_\Sigma$ indicate at which point the first order derivatives $A_k$ and $C_k$ of the functions $f$ and $h$ are computed. This notation is needed, if we want to use the extended Kalman in order to estimate how accurate our future state estimates will be without knowing the measurements yet. In such a situation, we can compute the iterates $\Sigma_k$ of the EKF by evaluating the recursion~\eqref{eq::ekf2} at predicted states $x_k$ in place of the state estimates $y_k$, which are in this case not available yet. This is possible, as the function $\mathcal F_{\Sigma}$ does not depend on $\eta_k$. In the following, we denote the true state sequence by $z_0, z_1, z_2, \ldots \in \mathbb R^{n_x}$. This sequence is defined by the recursion
\begin{eqnarray}
\forall k \in \mathbb N, \qquad z_{k+1} = f(z_k,u_k,w_k) \quad \text{with} \quad z_{0} = x_{0}^* \; .
\end{eqnarray}
Notice that this recursion depends on the disturbance sequence $w_{0}, w_1, w_2, \ldots$ as well as the true state $x_0^*$ of the system at time $0$, which are all unknown.

\bigskip
\begin{lemma}
\label{lem::EKF}
Let Assumptions~\ref{ass::differentiable},~\ref{ass::variances}, and~\ref{ass::boundedSupport} be satisfied, and let the sequences $y_0, y_1,y_2, \ldots$ and $\Sigma_0, \Sigma_1,\Sigma_2, \ldots$ be given by the extended Kalman filter recursions~\eqref{eq::ekf1} and~\eqref{eq::ekf2}, where the linearization points $x_k$ satisfy $x_k = z_k + \order{\gamma}$ while the sequence $u_k$ is given and exactly known. If the initial state estimate $y_0$ and the initial variance $\Sigma_0$ satisfy $\Sigma_0 = \order{\gamma^2}$ as well as
\[
\mathbb E \{ y_0 \} = x_0^* + \order{\gamma^2} \quad \text{and} \quad  \mathbb E \{ \left( y_0 - x_0^* \right)\left( y_0 - x_0^* \right)^\tr \} = \Sigma_0 + \order{\gamma^3} \; ,
\]
then the following statements hold on any finite time horizon, $k \in \{ 0, 1, \ldots, N \}$.

\bigskip
\begin{enumerate}

\item We have $\Sigma_k = \order{\gamma^2}$.

\item The expectation of the state estimate satisfies $\mathbb E \left\{ y_k \right\} = z_k + \order{\gamma^2}$.

\item The variance of the state estimate satisfies
\[
\mathrm{Var}(y_k - z_k) = \mathbb E \left\{ (y_k - z_k)(y_k - z_k)^\tr \right\} + \order{\gamma^3} = \Sigma_k + \order{\gamma^3} \; .
\]
\end{enumerate}
\end{lemma}

\bigskip
\proof
The statement of Lemma~\ref{lem::EKF} is---at least in very similar versions---well known in the literature~\cite{Stengel1994}. The main idea is to use induction. Firstly, all three statements of \change{Lemma~\ref{lem::EKF}} are true for $k = 0$. Next, under the assumption that these three statements are satisfied for a given $k$, they can be established for $k+1$: if $\Sigma_k = \order{\gamma^2}$ for a given $k$, then we also have
\[
\Sigma_{k+1} = F_{\Sigma}( x_{k}, u_{k}, \Sigma_{k} ) = \order{\Vert \Sigma_k \Vert_2} + \order{\Vert W \Vert_2} = \order{\gamma^2} \; ,
\]
i.e., the first statement holds for all $k \in \{ -N, \ldots, 0 \}$. The second and third statement follows from the fact the EKF computes the first and second order moments exactly, if the dynamic system $f$ and the measurement function $h$ are affine in $x$. Thus, we can apply the Taylor expansion techniques from Section~\ref{sec::prelim} to bound the contribution of higher order terms in order to show that the second and third statement also hold for $k+1$. The details of this argument can be found in~\cite{Haverbeke2011}.\hfill$\diamond$

\bigskip
\begin{Remark}
The statement of Lemma~\ref{lem::EKF} relies on the assumption that the horizon $N$ is finite, but the statement of this lemma is wrong in general for infinite horizons without further observability assumptions, as the EKF updates are divergent in general.\hfill$\diamond$
\end{Remark}

\bigskip
\begin{Remark}
Lemma~\ref{lem::EKF} can also be applied for the case that the state estimates and variances are computed with a moving horizon estimator (MHE)~\cite{Rawlings2009}, whose arrival cost is computed by extended Kalman filter updates. This leads to a refinement of the state estimates for nonlinear system. However, MHE with suitable arrival cost updates is in a first order approximation equivalent to EKF. A proof of this statement can be found in Chapter 3 of~\cite{Haverbeke2011}.\hfill$\diamond$
\end{Remark}
\bigskip

\section{Expected Loss of Optimality}
\label{sec::loss}

In Section~\ref{sec::motivating_example} we have discussed an example, where certainty equivalent MPC leads to a sub-optimal performance, as it tries to bring the system to a non-observable steady-state. The goal of this section is to analyze in detail how much performance is actually lost when running a certainty equivalent MPC to control a general nonlinear system in the presence of disturbances and inexact state estimates. Let $x_0^*$ denote the true (but unknown) value of the state at time $t = 0$. If we would know the future disturbance sequence $w^* = \change{ \left[ w_0^* \, w_1^* \, \ldots \, w_{N-1}^* \right] }$ in advance, we could compute the optimal input sequence $u^* = \change{ \left[ u_0^* \, u_1^* \, \ldots \, u_{N-1}^* \right] }$ as the solution of the optimal control problem
\begin{align}
\label{eq::mpcRep}
\begin{array}{rccl}
J^* = J_0(x_0^*, w_{[0]}^* ) &=& \underset{x,u}{\text{min}} \; & \sum_{k=0}^{N-1} l(x_k,u_k) + m(x_N) \\[0.3cm]
& & \text{s.t.} \; & \left\{
\begin{array}{l}
\forall k \in \{ 0, \ldots, N-1 \}, \\[0.16cm]
x_{k+1} = f(x_k,u_k) + w_k^* \\[0.16cm]
x_0 = x_0^* \; .
\end{array}
\right.
\end{array}
\end{align}
\change{In} practice, we neither know the true initial value $x_0^*$ nor can predict the future disturbance sequence $w^*$. Thus, the optimal input sequence $u^*$ cannot be computed. Recall that the control law, which is associated with the certainty equivalent MPC controller~\eqref{eq::mpc}, is given by
\[
\mu_k(\cdot) \; = \; \underset{u_k}{\mathrm{argmin}} \; G_k( \cdot, u_k, 0 ) \; ,
\]
where the functions $G_k$ have been defined \change{in~\eqref{eq::Gdef}} in Section~\ref{sec::mpc}. The true state of the closed loop MPC system at time $k$, denoted by $\xi_k^*(y^{[k-1]})$, is a function of the state estimate sequence
$$y^{[k-1]} = \change{ \left[ y_0 \, \ldots \, y_{k-1} \right] } \; .$$
It can be obtained recursively as the solution of the certainty equivalent MPC recursion
\begin{eqnarray}
\label{eq::trueMPC}
\begin{array}{rcl}
\xi_0^* &=& x_0^* \\[0.1cm]
\xi_1^*(y^{[0]}) &=& f(\xi_0^*,\mu_0(y_0))+w_0^* \\[0.1cm]
\xi_2^*(y^{[1]}) &=& f(\xi_1^*(y^{[0]}),\mu_1(y_1)) + w_1^* \\[0.1cm]
&\vdots& \\[0.1cm]
\xi_N^*(y^{[N-1]}) &=& f(\xi_{N-1}^*(y^{[N-2]}), \mu_{N-1}(y_{N-1})) + w_{N-1}^* \; .
\end{array}
\end{eqnarray}
In the following, we denote with
\begin{align}
\label{eq::subOptPerformance}
\begin{array}{l}
J_\text{cl}^*(y^{[N-1]}) \; = \; \sum_{k=0}^{N-1} l( \xi_k^*(y^{[k-1]}), \mu_k(y_k) ) + m(\xi_N^*(y^{[N-1]}))
\end{array}
\end{align}
the actual performance of the MPC controller. If uncertainties are present, the performance of the MPC controller can be expected to be worse than the best possible performance, given by
\begin{align}
\label{eq::OptPerformance}
\begin{array}{rcl}
J^* &=& \sum_{i=0}^{N-1} l( x_i^*, u_i^* ) + m(x_N^*) = J_\text{cl}^*( \, \change{ \left[ x_0^* \, x_1^* \, \ldots \, x_{N-1}^* \right] } \, ) \; .
\end{array}
\end{align}
We call the associated difference, given by
\[
\Delta^*(y) = J_\text{cl}^*(y) - J^* \, ,
\]
the ``loss of optimality''. In the following, we use the notation
\begin{eqnarray}
\label{eq::StageLoss}
\qquad \Delta_k^*(y^{[k]}) &=& G_k( \xi_k^*(y^{[k-1]}), \mu_k(y_k), w_{[k]}^* ) - G_k( \xi_k^*(y^{[k-1]}), \mu_k(\xi_k^*(y^{[k-1]})), w_{[k]}^* )
\end{eqnarray}
to denote the loss of optimality that is associated with the $k$-th suboptimal decision. The lemma below establishes an important relation between the function $\Delta^*$ and the function sequence $\Delta_0^*, \ldots, \Delta_{N-1}^*$.

\bigskip
\begin{lemma}
\label{lem::sum}
Let Assumption~\ref{ass::differentiable} be satisfied. We have $\Delta^*(y) = \sum_{k=0}^{N-1} \Delta_k^*(y^{[k]})$.
\end{lemma}

\bigskip
\proof
Recall that the functions $G_k$ are defined \change{by~\eqref{eq::Gdef}}. Together \change{with~\eqref{eq::trueMPC},} this definition implies
\begin{eqnarray}
\notag
G_k( \xi_k^*(y^{[k-1]}), \mu_k(y_k), w_{[k]}^* ) &=& J_{k+1}( f( \xi_{k}^*(y^{[k-1]}), \mu_k( y_k ) ) + w_{[k]}^*, w_{[k+1]}^* ) \\[0.16cm]
& & + l( \xi_k^*(y^{[k-1]}), \mu_k(y_k) ) \notag \\[0.16cm]
\label{eq::aux1}
&=& J_{k+1}( \xi_{k+1}^*(y^{[k]}), w_{[k+1]}^* ) + l( \xi_k^*(y^{[k-1]}), \mu_k(y_k) ) \; .
\end{eqnarray}
Moreover, since $\mu_k(\xi_k^*(y^{[k-1]}))$ is, by definition, a minimizer of the function $G_k( \xi_k^*(y^{[k-1]}), \cdot, w_{[k]}^* )$ the dynamic programming equation~\eqref{eq::dynProg} yields
\begin{eqnarray}
\label{eq::aux2}
G_k( \xi_k^*(y^{[k-1]}), \mu_k(\xi_k^*(y^{[k-1]})), w_{[k]}^* ) = J_{k}( \xi_{k}^*(y^{[k-1]}), w_{[k]}^* ) \; .
\end{eqnarray}
\change{Subtracting~\eqref{eq::aux1} and~\eqref{eq::aux2} yields}
\begin{eqnarray}
\Delta_k^*(y^{k}) &=& J_{k+1}( \xi_{k+1}^*(y^{[k]}), w_{[k+1]}^* ) - J_{k}( \xi_{k}^*(y^{[k-1]}), w_{[k]}^* ) + l( \xi_k^*(y^{[k-1]}), \mu_k(y_k) )\notag
\end{eqnarray}
for all $k \in \{ 0, 1, \ldots, N-1 \}$. Consequently, the sum over the terms $\Delta_k^*(y^{k})$ turns out to be a telescoping series, which can be simplified as follows:
\begin{eqnarray}
\sum_{k=0}^{N-1} \Delta_k^*(y^{[k]}) &=& J_{N}( \xi_{N}^*(y^{[N-1]}) ) - J_{0}( \xi_{0}^*, w_{[0]}^* ) + \sum_{k=0}^{N-1} l( \xi_k^*(y^{[k-1]}), \mu_k(y_k) ) \notag  \\[0.2cm]
&=& \sum_{k=0}^{N-1} l( \xi_k^*(y^{[k-1]}), \mu_k(y_k) ) + m(\xi_{N}^*(y^{[N-1]})) - J_{0}( x_{0}^*, w_{[0]}^* ) \notag \\[0.2cm]
&\overset{\eqref{eq::subOptPerformance}}{=}& J_\text{cl}^*(y) - J^* = \Delta^*(y) \; .
\end{eqnarray}
This is the statement of the lemma.\hfill$\diamond$

\bigskip
\noindent
The performance of an MPC controller that operates under inexact state measurements and process noise can never have a better performance than the utopian \change{controller~\eqref{eq::mpcRep}, which} knows the true current state and all future process disturbances. Thus, the loss of optimality function $\Delta^*$ must be non-negative. The following corollary summarizes this result and provides a formal proof.

\bigskip
\begin{corollary}
Let Assumption~\ref{ass::differentiable} be satisfied. We have $\Delta^*(y) \geq 0$ for all measurement sequences $y = y^{[N-1]}$. 
\end{corollary}

\bigskip
\proof
Definition~\eqref{eq::StageLoss} implies that the functions $\Delta_k$ are non-negative, as $\mu_k(\xi_k^*(y^{[k-1]}))$ is a minimizer of the function $G_k( \xi_k^*(y^{[k-1]}), \cdot, w_{[k]}^* )$. Thus Lemma~\ref{lem::sum} implies that the function $\Delta^*$ is non-negative, too.
\hfill$\diamond$

\bigskip
\noindent
In the following, we focus on the derivation of a tractable second order approximation of the expected loss of optimality, given by
\[
\mathbb E \left\{ \Delta^*(y) \right\} = \sum_{k=0}^{N-1} \mathbb E \left\{ \Delta_k^*(y^{[k]}) \right\} \; .
\]
Here, the expected value is computed under the assumption that the state estimates $y$ are computed by the EKF (or locally equivalent MHE) that is described in Section~\ref{sec::mhe}. On the first view, we might expect that the computation of such a second order expansion is rather cumbersome and computationally expensive, as the state estimates $y_k$ are in general correlated in time. Nevertheless, the following theorem provides a surprisingly simple recursion technique for computing a second order approximation of the expected loss of optimality requiring $\order{N}$ operations only. It is based on the following notation.

\bigskip
\begin{itemize}\itemsep4pt

\item The matrices $\Omega_k = [p_k,P_k]$ denote the solution of the backward Riccati recursion
\begin{eqnarray}
\label{eq::backaux}
\Omega_k &=& \mathcal F_{\Omega}(x_k,u_k,\Omega_{k+1}) \quad \text{with} \quad \Omega_N = \mathcal M(x_N)
\end{eqnarray}
for $k \in \{ 0, \ldots, N-1 \}$, which has been introduced in Section~\ref{sec::mpc}.

\item We use the shorthand notation
\begin{eqnarray}
\Phi_k &=& G_k^{xu} \left( G^{uu}_k \right)^{-1} G^{ux}_k \notag \\[0.16cm]
&=& \left( S_k + B_k^\tr P_{k+1} A_k + L_k \cdot p_{k+1} \right)^\tr \cdot \left( R_k + B_k^\tr P_{k+1} B_k + M_k \cdot p_{k+1} \right)^{-1} \notag \\[0.16cm]
& & \cdot \left( S_k + B_k^\tr P_{k+1} A_k + L_k \cdot p_{k+1} \right) . \notag
\end{eqnarray}
Notice that the matrices $\Phi_k$ need to be computed anyhow during the evaluation of the backward recursion~\eqref{eq::backaux}, as the recursion for the matrices $P_k$ can be written in the explicit form
\begin{eqnarray}
P_k = \underbrace{\left( Q_k + A_k^\tr P_{k+1} A_k + K_k \cdot p_{k+1} \right)}_{G_k^{xx}} - \Phi_k \; .
\end{eqnarray}

\item Finally, the matrices $\Sigma_k$ are given by the forward Riccati recursion
\begin{eqnarray}
\Sigma_{k+1} &=& \mathcal F_{\Sigma}( x_{k}, u_{k}, \Sigma_{k} ) \quad \text{with} \quad \Sigma_0 = \hat \Sigma_0
\end{eqnarray}
where $\hat \Sigma_0 = \mathbb E \{ (y_0 - x_0^*)(y_0 - x_0^*)^\tr \} + \order{\gamma^3}$ denotes an approximation of the variance of the initial state estimate $y_0$ as introduced in Section~\ref{sec::mhe}.
\end{itemize}

\bigskip
\noindent
Theorem~\ref{thm::expansion} in combination with Corollary~\ref{cor::expansion} are one of the core technical contributions of this paper.

\bigskip
\begin{theorem}
\label{thm::expansion}
Let Assumptions~\ref{ass::differentiable},~\ref{ass::variances},~\ref{ass::boundedSupport}, and~\ref{as::regularity} be satisfied, then we have
$$\mathbb E \, \left\{ \Delta_k^*(y^{[k]}) \right\} \; = \; \frac{1}{2} \mathrm{Tr}( \Phi_k \Sigma_k ) + \mathbf{O}\left( \gamma^3 \right)$$
for all $k \in \{ 0, \ldots, N-1 \}$. Here, the matrices $\Omega_k, \Phi_k$, and $\Sigma_k$ are computed as outlined above and under the assumption that the linearization point $(x_k,u_k)$, at which the forward and backward Riccati recursions for the sequences $\Omega_k$ and $\Sigma_k$ as well as the matrices $\Phi_k$ are evaluated, satisfies $(x_k,u_k) = (y_k,\mu_k(y_k)) + \order{\gamma}$.
\end{theorem}

\bigskip
\proof
In order to understand the following proof it is helpful to recall that the state estimates $y_1, \ldots, y_k$ are correlated random variables. Thus, if we compute the expectation with respect to the most recent state estimate $y_k$ first,
\begin{align}
\begin{array}{rcl}
\mathbb E \left\{ \Delta_k^*(y^{[k]}) \right\} &=& \underset{y^{[k-1]}}{\mathbb E} \left\{ \underset{y_k}{\mathbb E} \left\{ \Delta_k^*(y^{[k]}) \right\} \right\} \; ,
\end{array}
\end{align}
we should keep in mind that $\underset{y_k}{\mathbb E} \left\{ \Delta_k^*(y^{[k]}) \right\}$ denotes the conditional expectation value of $\Delta_k^*(y^{[k]})$ over the variable $y_k$, which does depend on the values of the previous state estimates
$y^{[k-1]} = \left( y_0,\ldots,y_{k-1} \right) \; .$
Since $\mu_k$ is an optimizing control law, the first order necessary optimality condition
\[
\left. \frac{\partial G_k(\xi_k^*(y^{[k-1]}),u,w_{[k]}^*)}{\partial u} \right|_{u = \mu_k(\xi_k^*(y^{[k-1]}))} = 0
\]
must be satisfied. Thus, it follows \change{from~\eqref{eq::StageLoss}} that the equations
\begin{align}
0 = \left. \Delta_k^*(y^{[k]}) \right|_{y_k = \xi_k^*(y^{[k-1]})}  \quad \text{and} \quad 0 &=& \left. \frac{\partial}{\partial y_k} \Delta_k^*(y^{[k]}) \right|_{y_k = \xi_k^*(y^{[k-1]})}
\end{align}
hold independently of the choice of $y^{[k-1]}$. Next, we introduce the notation
$$\Gamma_k(y^{[k-1]}) = \left. \nabla_{y_k}^2 \Delta_k^*(y^{[k]}) \right|_{y_k = \xi_k^*(y^{[k-1]})}$$
to denote the Hessian matrix of the function $\Delta_k^*$. This Hessian matrix can be worked out by applying the implicit function theorem in combination with the chain rule,
\begin{eqnarray}
\begin{array}{l}
\Gamma_k(y^{[k-1]}) = \left. \dfrac{\partial^2 G_k(x,u,w^*)}{\partial x \partial u} \right|_{\substack{
u = \mu_k(\xi_k^*(y^{[k-1]}))\\
x = \xi_k^*(y^{[k-1]}) \; \; \quad
}} \cdot \left( \left. \dfrac{\partial^2 G_k(x,u,w^*)}{\partial u \partial u} \right|_{\substack{
u = \mu_k(\xi_k^*(y^{[k-1]}))\\
x = \xi_k^*(y^{[k-1]}) \; \; \quad
}} \right)^{-1} \notag \\[1.0cm]
\hspace{2.5cm} \cdot \left. \dfrac{\partial^2 G_k(x,u,w^*)}{\partial u \partial x} \right|_{\substack{
u = \mu_k(\xi_k^*(y^{[k-1]}))\\
x = \xi_k^*(y^{[k-1]}) \; \; \quad
}} \; . \notag
\end{array}
\end{eqnarray}
Fortunately, this equation can be written in the form
\[
\Gamma_k(y_0,\ldots,y_{k-1}) = \Phi_k + \order{\gamma}
\]
with $\Phi_k = G_{k}^{xu} \left( G_{k}^{uu} \right)^{-1} G_{k}^{ux}$, since we assume that the evaluation point $(x_k,u_k)$ satisfies
$$(x_k,u_k) = \left( y_k,\mu_k(y_k) \right) + \order{\gamma} \; .$$
Moreover, it follows from the consideration in Section~\ref{sec::mhe} that we have
$$\mathbb E \left\{ \left( y_k - \xi_k^*(y^{[k-1]} \right)\left( y_k - \xi_k^*(y^{[k-1]} \right)^\tr \right\} = \Sigma_k + \order{\gamma^3} \; .$$
Consequently, by using the second order moment expansion technique from Section~\ref{sec::prelim}, the term $\underset{y_k}{\mathbb E} \left\{ \Delta_k^*(y^{[k]}) \right\}$ can be written in the form
\[
\underset{y_k}{\mathbb E} \left\{ \Delta_k^*(y^{[k]}) \right\} = \frac{1}{2} \, \mathrm{Tr} \left( \Phi_k \Sigma_k \right) + \order{\gamma^3} \; .
\]
This equation is fortunate, since it implies that the conditional expectation $\underset{y_k}{\mathbb E} \left\{ \Delta_k^*(y^{[k]}) \right\}$ is in a second order approximation independent of $y^{[k-1]}$. Thus, the expectation value over all random variables can be evaluated as
\begin{align}
\begin{array}{rcl}
\mathbb E \left\{ \Delta_k^*(y^{[k]}) \right\} \; = \; \underset{y^{[k-1]}}{\mathbb E} \left\{ \underset{y_k}{\mathbb E} \left\{ \Delta_k^*(y^{[k]}) \right\} \right\} \; = \; \frac{1}{2} \, \mathrm{Tr} \left( \Phi_k \Sigma_k \right) + \order{\gamma^3} \; . \notag
\end{array}
\end{align}
This completes the proof of Theorem~\ref{thm::expansion}.\hfill$\diamond$

\bigskip
\begin{corollary}
\label{cor::expansion}
If the conditions from Theorem~\ref{thm::expansion} are satisfied, then the expected loss of optimality is given by
\[
\mathbb E \{ \Delta^*(y) \} \; = \; \frac{1}{2} \sum_{k=0}^{N-1} \, \mathrm{Tr}( \Phi_k \Sigma_k ) + \mathbf{O}\left( \gamma^3 \right) \; .
\]
\end{corollary}

\proof
The statement of this corollary is a trivial consequence of Theorem~\ref{thm::expansion} and Lemma~\ref{lem::sum}.\hfill$\diamond$

\section{Self-Reflective Model Predictive Control}
\label{sec::rhc}

This section proposes a receding horizon controller, whose objective is to minimize the sum of the certainty equivalent control performance and a second order approximation of the associated expected loss of optimality due to inaccurate state estimates. For this aim, we recall that the second order approximation $\Sigma_k$ of the variance of the future state estimates satisfies the forward recursion
\[
\Sigma_{k+1} = \mathcal F_{\Sigma}( x_k, u_k, \Sigma_k ) \quad \text{with} \quad \Sigma_0 = \hat \Sigma_0
\]
for all $k \in \{ 0, \ldots, N-1 \}$. The associated initial value is assumed to satisfy
$$\hat \Sigma_0 = \mathbb E \{ (y_0 - x_0^*)(y_0 - x_0^*)^\tr \} + \order{\gamma^3} \; ,$$
where $y_0$ denotes the current state estimate and $x_0^*$ the true state at time $t = 0$. As in the previous section $\Omega_k = [p_k,P_k]$ denotes the state of the backward Riccati recursion~\eqref{eq::forwardpP} with the parametric final state $\Omega_N = \mathcal M(x_N)$. Motivated by Theorem~\ref{thm::expansion}, we introduce the notation
\[
\mathcal L(x_k,u_k,\Sigma_k,\Omega_{k+1}) \; = \; \frac{1}{2} \mathrm{Tr}( \Phi_k \Sigma_k )
\]
to denote the second order approximation of the expected loss of optimality that is associated with the $k$-th future decision. Here, the matrix $\Phi_k$ is a function of the linearization point $(x_k,u_k)$ and the state $\Omega_{k+1}$ of the backward Riccati recursion. The focus of this section is on the parametric optimal control problem
\begin{eqnarray}
\begin{array}{rccl}
J_i^\mathrm{SR}(x_i,\Sigma_i) &=& \underset{x,u,\Sigma,\Omega}{\text{min}} \; & \sum_{k=i}^{N-1} \left\{ l(x_k,u_k) + \change{\alpha} \mathcal L(x_k,u_{k},\Sigma_k,\Omega_{k+1}) \right\} + m(x_N) \\[0.3cm]
& & \text{s.t.} \; &
\left\{
\begin{array}{l}
\forall k \in \{ i, \ldots, N-1 \},\\[0.16cm]
\begin{array}{rcl}
x_{k+1} &=& f(x_k,u_k), \\[0.16cm]
\Sigma_{k+1} &=& \mathcal F_{\Sigma}( x_k, u_k, \Sigma_k ), \\[0.16cm]
\Omega_{k}      &=& \mathcal F_{\Omega}( x_k, u_{k}, \Omega_{k+1} ), \\[0.16cm]
\Omega_N &=& \mathcal M(x_N) \, .
\end{array}
\end{array}
\right.
\end{array}
\label{eq::srmpc}
\end{eqnarray}
\change{Here, $\alpha \geq 0$ is a non-negative tuning parameter that can be used to trade-off nominal performance versus expected loss of optimality.}
The full self-reflective MPC loop can be outlined as follows.

\bigskip
\begin{enumerate}
\item Wait for the measurement $\eta_0$.
\item Update the current state estimate and variance by using an EKF
\begin{eqnarray}
\notag \hat y &\leftarrow& \mathcal F_y( \hat y, u_{-1}, \hat y, \eta_0 ) \; ,\\[0.16cm]
\notag \hat \Sigma &\leftarrow& \mathcal F_{\Sigma}( \hat y, u_{-1}, \hat \Sigma ) \; .
\end{eqnarray}
\item Solve Problem~\eqref{eq::srmpc} for $i=0$, $x_0 = \hat y$, and $\Sigma_0 = \hat \Sigma$.
\item Send the solution for $u_0$ to the real process.
\item Reset the time index and continue with Step~1. 
\end{enumerate}

\bigskip
\noindent
In the following, we denote with $\mu_i^\mathrm{SR}(x_i,\Sigma_i)$ the optimal solution of Problem~\eqref{eq::srmpc} for the optimization variable $u_i$ in dependence on $x_i$ and $\Sigma_i$. The true closed-loop trajectory is then given by the coupled recursion
\begin{eqnarray}
\zeta_{k+1}^*( y^{[k]} ) &=& f( \zeta_{k}^*( y^{[k-1]} ), \mu_k^\mathrm{SR}(y_k,\Sigma_k^*(y^{[k-1]}))  ) + w_k^* \quad \text{with} \quad \zeta_0^* = x_0^* \\[0.16cm]
\Sigma_{k+1}^*(y^{[k]}) &=& F_\Sigma( y_k, \mu_k^\mathrm{SR}(y_k,\Sigma_k^*(y^{[k-1]})), \Sigma_k^*(y^{[k-1]}) ) \quad \text{with} \quad \Sigma_0^* = \Sigma_0^* \; .
\end{eqnarray}
The loss of optimality of the above outlined controller is given by
\[
\Delta_{\mathrm{SR}}^*(y^{[N-1]}) = \sum_{k=0}^N l( \zeta_{k}^*( y^{[k-1]} ), \mu_k^\mathrm{SR}(y_k,\Sigma_k^*(y^{[k-1]})) ) + m(\zeta_{N}^*( y^{[N-1]} )) - J^* \; .
\]
Finally, the expected loss of optimality is denoted by $\mathbb E_y \left\{ \Delta_{\mathrm{SR}}^*(y^{[N-1]}) \right\}$. The following corollary is a simple consequence of Theorem~\eqref{thm::expansion}.


\bigskip
\begin{corollary}
Let Assumptions~\ref{ass::differentiable},~\ref{ass::variances},~\ref{ass::boundedSupport}, and~\ref{as::regularity} be satisfied, then the above outlined controller from Problem~\eqref{eq::srmpc} is (approximately) self-reflective in the sense that it optimizes the sum of its nominal performance and a second order approximation of its own expected loss of optimality. In other words, if $(x,u,\Sigma,\Omega)$ denotes an optimal solution of Problem~\eqref{eq::srmpc} for $i=0$ and
$$\mathbb E \{ \hat y \} = x_0^* + \order{\gamma^2} \quad \text{and} \quad \hat \Sigma = \mathbb E \{ (y_0 - x_0^*)(y_0 - x_0^*)^\tr \} + \order{\gamma^3} \; ,$$
then we have
\[
\mathbb E_y \left\{ \Delta_{\mathrm{SR}}^*(y^{[N-1]}) \right\} = \sum_{k=0}^{N-1} \mathcal L(x_k,u_{k},\Sigma_k,\Omega_{k+1}) + \order{\gamma^3}\; .
\]
\end{corollary}

\bigskip
\proof
Notice that Problem~\eqref{eq::srmpc} and Problem~\eqref{eq::mpc} differ for small $\gamma$ only by the small (but not insignificantly small) additional term $\mathcal L(x_k,u_{k},\Sigma_k,\Omega_{k+1})$, which satisfies
\[
\change{\alpha} \mathcal L(x_k,u_{k},\Sigma_k,\Omega_{k+1}) = \order{\gamma^2} \; .
\]
Thus, due to the regularity requirement from Assumption~\ref{as::regularity}, this additional term perturbs the control law by terms of order $\order{\gamma^2}$ only,
\[
\mu_k^{\mathrm{SR}}(x,\Sigma) = \mu_k(x) + \order{\gamma^2}
\]
as long as $\Sigma = \order{\gamma^{2}}$. An analogous conclusion holds for the state trajectory, which implies
\begin{eqnarray}
\mathbb E_y \left\{ \Delta_{\mathrm{SR}}^*(y^{[N-1]}) \right\} = \mathbb E_y \left\{ \Delta^*(y^{[N-1]}) \right\} + \order{\gamma^4} \; .
\end{eqnarray}
We know from Theorem~\ref{thm::expansion} that the expected loss of optimality of the certainty equivalent MPC is given by $\mathbb E_y \left\{ \Delta^*(y^{[N-1]}) \right\} = \mathcal L(x_k,u_{k},\Sigma_k,\Omega_{k+1}) + \order{\gamma^3}$ and that this equation holds independently of the linearization point $(x_k,u_k)$ as long as we do not move out of a $\order{\gamma}$ neighborhood of the true optimal state and input sequence. This proves the statement of the corollary.\hfill$\diamond$

\bigskip
\begin{Remark}
Problem~\eqref{eq::srmpc} is not a standard optimal control problem in the sense that the stage cost contains a term that couples the $k$-th state $x_k,\Sigma_k$ and the $(k+1)$-th state $\Omega_{k+1}$. However, Problem~\eqref{eq::srmpc} is still band-structured and there exist structure-exploiting solvers for implicit dynamic systems~\cite{Biegler1991}, which can be used to solve problems of the form~\eqref{eq::srmpc} with run-time complexity $\order{N}$.\hfill$\diamond$
\end{Remark}

\bigskip
\begin{Remark}
As already mentioned in the introduction, existing methods on persistently exciting MPC controllers such as~\cite{Hovd2005} are similar to the proposed controller in the sense that they also add a term that minimizes the trace or other scalar measures of the predicted variance state $\Sigma_k$ of future state estimates. However, the main contribution of this paper is that we compute weighting matrices $\Phi_k$ such that the penalty terms can be interpreted as the expected loss of optimality of the controller. Notice that not only the variances $\Sigma_k$ of future state estimates but also the weighting matrices $\Phi_k$ may depend on the control input $u$. If we would only add a term of the form $\mathrm{Tr}(\Sigma_k)$ to the objective, this might lead to more accurate state estimates but it does not necessarily improve the expected control performance, if the associated control decisions lead to an increase of the eigenvalues of the matrix $\Phi_k$.\hfill$\diamond$
\end{Remark}

\bigskip
\begin{Remark}
The above analysis can be generalized for the case that control constraints are present. One obvious solution for including such constraints is to relax them and add barrier functions~\cite{Nocedal2006} to the Lagrange term $l$, which penalize constraint violations. If such ``soft-constraints'' are acceptable, the barrier functions can be chosen to be smooth such that the above analysis can be applied. Otherwise, if Problem~\eqref{eq::mpc} is augmented with hard constraints on the control input, the cost-to-go functions $J_i$ and, consequently, also the associated functions $G_i$ are only piecewise differentiable. Thus, second order expansions of the expected loss of optimality are only possible in the neighborhood of regular points, where no active set changes are expected. If we are at a non-differentiable point of the cost-to-go function, where active set changes have to be taken into account, estimates of the expected loss of optimality can be evaluated by using Monte-Carlo simulations, which are, however, much more expensive than Taylor expansion based second order estimates. Other strategies for dual control in the presence of control constraints may be possible and an interesting research direction, but are beyond the scope of the current paper.
\end{Remark}

\bigskip
\noindent
A rigorous analysis of the limit behavior and the associated question about the stability of the proposed self-reflective controller~\eqref{eq::srmpc} is an important open problem, which is, however, beyond the scope of this paper. Notice that a rigorous stability analysis of existing persistently exciting MPC, let alone dual controllers, is at the current status of research equally unavailable. However, there exists a number of articles that address the stability of general economic MPC controllers with and without terminal costs~\cite{Amrit2011,Angeli2012,Diehl2011,Gruene2011,Houska2015a,Zanon2013}. Also notice that the presented self-reflective MPC scheme does not address the question how the uncertainty affects the states, which is in contrast to  tube based robust MPC~\cite{Rakovic2012}. An analysis of the robustness properties of self-reflective MPC might be an interesting topic for future research.

\section{Case Study}
\label{sec::numerics}

This section presents a numerical case study for a predator-prey system of the form
\change{
\begin{eqnarray}
\dot z_1(t) &=& z_1(t) - z_1(t) z_2(t) + w_1(t) \; , \\[0.1cm]
\dot z_2(t) &=& - z_2(t) + z_1(t) z_2(t) + u(t) z_3(t) z_2(t) + w_2(t) \; , \\[0.1cm]
\dot z_3(t) &=& - z_3(t) + 0.5 + w_3(t) \; .
\end{eqnarray}
}
This model is a slightly modified version of the predator-prey optimal experiment design benchmark problem from~\cite{Sager2013}. Here, the states $z_1$ and $z_2$ can be interpreted as the population densities of a prey and a predator species, respectively. We can use the control input $u(t)$ to feed the predator with random success rate $z_3(t) \in [0,1]$. Notice that the approximately Brownian noise $w_3(t) \in [-0.5,0.5]$ is assumed to have bounded support such that $z_3(t)$ is for all $t > 0$ in $[0,1]$ as long as $z_3(0) \in [0,1]$. The functions $w_1(t)$ and $w_2(t)$ model additional external approximately Brownian noise affecting the population density of the predator and prey, respectively. \change{We} assume that we can only measure the population of the prey,
\[
h(z(t)) = z_1(t) \; ,
\]
but we have neither direct measurements of $z_2$ nor of $z_3$. Now, \change{the} nominal objective has the form
\[
\int_0^T \left( (z_1(t)-1)^2 + (z_2(t)-1)^2 + u(t)^2 \right) \, \mathrm{d}t + \underbrace{(z(T)-z_\mathrm{s})^\tr P_T (z(T)-z_\mathrm{s})}_{= m(z(T))} \; ,
\]
i.e., we would like to keep \change{the state $z = [ z_1 \, z_2 \, z_3]^\tr$ close to the steady state $z_\mathrm{s} = [1 \, 1 \, 0.5 ]^\tr$}. Notice that the system is not observable at the steady state, as the corresponding steady-state control input is given by $u_\text{s}=0$. The weighting \change{matrix $P_T \in \mathbb S_{+}^3$}
is constant and chosen in such a way that the terminal cost $m$ can be interpreted as an \change{approximation of the infinite horizon cost~\cite{Rawlings2009}.} In this example, we choose the prediction horizon $T = 10$. \change{Moreover, we use piecewise constant control and uncertainty discretization in combination with Euler's integration method such that the above differential equation can be approximated by a discrete-time system, which is affine in $w$. The sampling time of the discretization scheme and the controller are set to $\delta = 10^{-2}$ such that the discrete-time horizon of the associated self-reflective model predictive control problem~\eqref{eq::srmpc} is given by $N = T/\delta = 10^3$. One might argue that this discretization is neither efficient nor accurate. However, the focus of this paper is on discrete-time systems and therefore an analysis of the discretization error or more efficient discretization schemes is beyond the scope of this paper.} The matrices $V$ and $W$ are in our case study given by
\[
V = \change{\frac{1}{\delta}} \quad \text{and} \quad W = \frac{1}{\delta} \; \mathrm{diag}( \, \change{ 0.01, \,  0.01, \, 0.025 } \, ) \; .
\]
Moreover, the closed-loop system is started with
\change{ 
\[
\hat y_0 \, = \, \left(
\begin{array}{c}
1 \\
1 \\
0.5
\end{array}
\right) \qquad \text{and} \qquad
\hat \Sigma_0 \, = \, \mathrm{diag}( \, 0.1, \, 0.1, \, 0.1 \, ) \; .
\]
}
We \change{first} simulate the self-reflective MPC closed-loop system without applying process noise or measurement errors (but the controller does not know this).
\begin{figure}[t]
\begin{center}
\includegraphics[scale=0.65]{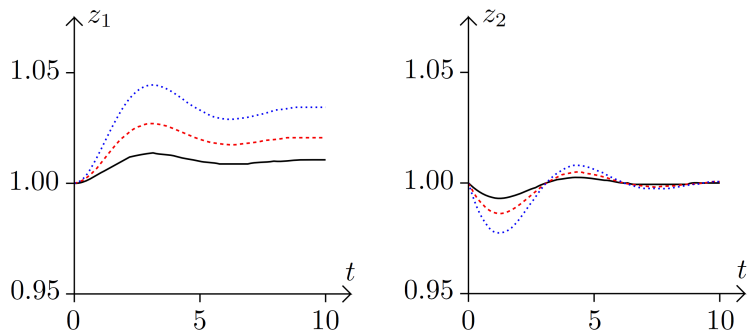}
\end{center}
\caption{\change{The closed-loop state trajectories $z_1$ and $z_2$ for $\alpha = \frac{1}{2}$ (black solid line), $\alpha = 1$ (red dashed line), and $\alpha = 2$ (blue dotted line) without noise.}}
\label{fig::1}
\end{figure}
Figure~\ref{fig::1} shows the numerical solution for the closed-loop trajectories for the predator and prey population densities \change{in dependence} on time (iteration index) \change{for different trade-off parameters $\alpha$. The numerical solution has been found by using the software ACADO Toolkit~\cite{Houska2011}.} Notice that the self-reflective MPC controller does not attempt to track the desired steady state \change{$z_\mathrm{s}$}. \change{The closed-loop state trajectory for the prey population converges to a value, which is larger than $1$. Moreover, if we increase the tuning parameter $\alpha$, the excitation of the system states increases. The table below lists the associated numerical values for the expected loss of optimality in dependence on $\alpha$ evaluated at the predicted trajectory.}

\begin{center}
\change{
\begin{tabular}{|c|c|}
\hline
$\alpha$ & $\sum_{k=0}^{N-1} \mathcal L(x_k,u_{k},\Sigma_k,\Omega_{k+1})$ \\[0.1cm]
\hline
$\frac{1}{2}$ & $\approx 2.94$ \\[0.1cm]
$1$ & $\approx 2.84$ \\[0.1cm]
$2$ & $\approx 2.80$ \\
\hline
\end{tabular}
}
\end{center}

\bigskip
\change{Clearly, the expected loss of optimality decreases for larger $\alpha$.}
\begin{figure}[t]
\begin{center}
\includegraphics[scale=0.65]{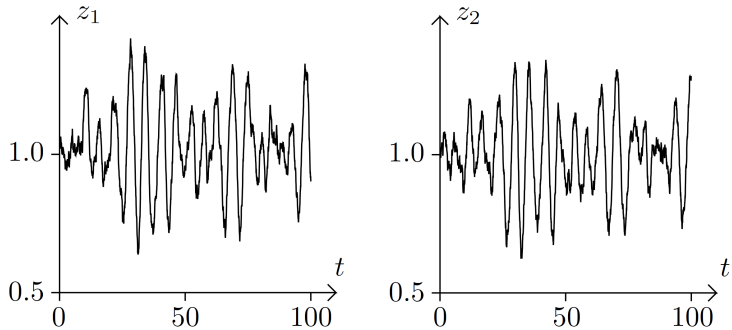}
\end{center}
\caption{\change{The closed-loop state trajectories $z_1$ and $z_2$ in the presence of random measurement errors and random process noise.}}
\label{fig::2}
\end{figure}
\change{Figure~\ref{fig::2} shows a long-term simulation of the closed-loop state trajectories for $\alpha = 1$ in the presence of random measurement errors and random process noise. Notice that the self-reflective MPC controller leads to reasonable control performance, which is in contrast to the closed-loop response of nominal MPC, which turns out to be divergent in this case study.}

\section{Conclusions}
\label{sec::conclusions}

In this paper a self-reflective model predictive controller has been proposed. Important contributions are Lemma~\ref{lem::sum} and Theorem~\ref{thm::expansion} establishing a recursion scheme for computing tractable second order expansions of the expected loss of optimality that is associated with certainty equivalent MPC schemes. This recursion scheme has been used to develop a model based controller, which is self-reflective in the sense that it minimizes the sum of its nominal performance and a second order approximation of its own expected loss of optimality due to inexact future measurements. Moreover, a small but non-trivial case study has illustrated that the controller is capable of exciting the system in order to improve future state estimate while maintaining economic control performance.

\section{Acknowledgments}
\label{sec::acknows}
This research was supported by National Natural Science Foundation China (NSFC), Nr.~61473185, ShanghaiTech University, Grant-Nr.~\mbox{F-0203-14-012}, as well as PFV/10/002 (OPTEC), FWO KAN2013 1.5.189.13,
FWO-G.0930.13 and BelSPO: IAP VII/19 (DYSCO). DT thanks KU Leuven for a postdoctoral grant.

\end{document}